\documentclass[11pt]{article}
\usepackage{amssymb,amsfonts,amsmath,latexsym,epsf,tikz,url}

\newtheorem{theorem}{Theorem}[section]

\newtheorem{corollary}[theorem]{Corollary}
\newtheorem{lemma}[theorem]{Lemma}
\newtheorem{remark}[theorem]{Remark}
\newtheorem{example}[theorem]{Example}
\newcommand{\proof}{\noindent{\bf Proof.\ }}
\newcommand{\qed}{\hfill $\square$\medskip}

\begin{document}

\title{Distinguishing number and distinguishing index of join of two graphs}

\author{
Saeid Alikhani  $^{}$\footnote{Corresponding author}
\and
Samaneh Soltani
}

\date{\today}

\maketitle

\begin{center}
Department of Mathematics, Yazd University, 89195-741, Yazd, Iran\\
{\tt alikhani@yazd.ac.ir, s.soltani1979@gmail.com}
\end{center}

%%%%%%%%%%%%%%ABSTRACT%%%%%%%%%%%%%%%%%%%%%%%%%%%%%%%%%%%%%%%%%%%%%%%%%%%%%%%%%%%%
\begin{abstract}
The distinguishing number (index) $D(G)$ ($D'(G)$) of a graph $G$ is the least integer $d$
such that $G$ has an vertex labeling (edge labeling)  with $d$ labels  that is preserved only by a trivial
automorphism. In this paper  we study the distinguishing number and the distinguishing index of join of two graphs $G$ and $H$, i.e., $G+H$. 
We
 prove that $0\leq D(G+H)-max\{D(G),D(H)\}\leq z$, where $z$ is depends of the number of some induced subgraphs generated by some suitable partitions of $V(G)$ and $V(H)$.  Also, we prove that if $G$ is a connected graph of order $n \geq 2$, then $D'(G+ \cdots +G)=2$, except $D'(K_2+K_2)=3$.
 \end{abstract}


\medskip
\noindent{\bf AMS Subj.\ Class.:} 05C15, 05E18

%%%%%%%%%%%%%%%%%%%%%%%%%%%%%%%%%%%%%%%%%%%%%%%%%%%%%%%%%%%%%%%%%%%%%%%%%%%%%%%%%
%%%%%%%%%%%%%%%%%%%%%%%%%%%%%%%%%%%%%%%%%%%%%%%%%%%%%%%%%%%%%%%%%%%%%%%%%%%%%%%%%
\section{Introduction}

Let $G = (V ,E)$ be a graph with $n$ vertices.  We use the standard graph notation (\cite{Sandi}). In particular, ${\rm Aut}(G)$ denotes the automorphism group of $G$.  
A labeling of $G$, $\phi : V \rightarrow \{1, 2, \ldots , r\}$, is  $r$-distinguishing, 
if no non-trivial  automorphism of $G$ preserves all of the vertex labels.
Formally, $\phi$ is $r$-distinguishing if for every non-trivial $\sigma \in {\rm Aut}(G)$, there
exists $x$ in $V$ such that $\phi(x) \neq \phi(\sigma x)$. 
The distinguishing number of a graph $G$ is the minimum number $r$ such that $G$ has a labeling that is $r$-distinguishing.  
This number was defined by Albertson and Collins \cite{Albert}. Similar to this definition, Kalinkowski and Pil\'sniak \cite{Kali1} have defined the distinguishing index $D'(G)$ of $G$ which is  the least integer $d$
such that $G$ has an edge colouring   with $d$ colours that is preserved only by a trivial
{\rm Aut}omorphism. Observe that  $D(G) = 1$ for the asymmetric graphs $G$ and  $D(G) = \vert V(G) \vert$,  if and only if $G = K_n$. It is immediate that $D(P_n) = 2$ for $n\geq 2$, where $P_n$ is the $n$-vertex path. A classical result gives that for the cycle with $n$ vertices, $C_n$, $D(C_n) = 3$ if $n = 3,4,5$ and $D(C_n) = 2$ for $n\geq 6$. Also for complete bipartite graph when $q > p$, $D(K_{p,q}) = q$,  $D(K_{n,n}) = n + 1$ for $n\geq 3$,  for the $n$-cube $Q_n$, $D(Q_n)=2$, for $n\geq 4$ and $D(Q_n)=3$ for $n=2,3$ (\cite{bogs}). 
The distinguishing index of some graphs was exhibited in \cite{soltani,Kali1}.
The distinguishing number and index  of the Cartesian product and the Cartesian  powers of graphs has been
thoroughly investigated (\cite{Albert2005,Klavzar,fish}). 
Pil\'sniak studied the Nordhaus-Gaddum bounds for the distinguishing index in \cite{nord}. Also the distinguishing number of the hypercube has been investigated 
in \cite{bogs}.  
Recently, we studied the distinguishing number and distinguishing index of corona product  of two graphs (\cite{soltani}). 

\medskip
 We say that $G= (V,E)$ is a join graph if $G$ is the complete union of two graphs $G_1 = (V_1,E_1)$ and $G_2 = (V_2,E_2)$. In other words, $V = V_1 \cup V_2$ and $E = E_1 \cup E_2 \cup \{uv \vert u \in V_1, v \in V_2\}$. If $G$ is
the join graph of $G_1$ and $G_2$, we write $G = G_1 +G_2$. 
For simple connected  graph $G$, and $v\in V$, the neighborhood of $v$ is the set $N_G(v) = N(v) = \{u \in V (G) : uv \in E(G)\}$. The nonadjacent vertices to $v$ in $G$ is $ V(G)\setminus N(v)$ and denoted by $\overline{N(v)}$. A subgraph $H$ of $G$ is an induced subgraph if two vertices of $V(H)$ are adjacent in $H$ if and only if they are adjacent in $G$. We denote the induced subgraph by a set $X\subseteq V$, by $G[X]$.

In the next section, we study the  distinguishing number  of join of two graphs. In Section 3, we present two upper bounds for the  distinguishing index of  the join of two graphs and show that they are sharp.    

%%%%%%%%%%%%%%%%%%%%%%%%%%%%%%%%%%%%%%%%%%%%%%%%%%%%%%%%%%%%%%%%%%%%%%%%%%%%%%%%%
%%%%%%%%%%%%%%%%%%%%%%%%%%%%%%%%%%%%%%%%%%%%%%%%%%%%%%%%%%%%%%%%%%%%%%%%%%%%%%%%%
\section{The distinguishing number of the join of two graphs}
%%%%%%%%%%%%%%%%%%%%%%%%%%%%%%%%%%%%%%%%%%%%%%%%%%%%%%%%%%%%%%%%%%%%%%%%%%%%%%%%%
%%%%%%%%%%%%%%%%%%%%%%%%%%%%%%%%%%%%%%%%%%%%%%%%%%%%%%%%%%%%%%%%%%%%%%%%%%%%%%%%%
In this section,  we study the  distinguishing number of join of two graphs. We begin  with the following theorem which gives
 a lower bound for the distinguishing number of join of two graphs: 
\begin{theorem}\label{thh5}
Let $G_1$ and $G_2$ be two  connected graphs. Then 
$${\rm max}\{D(G_1),D(G_2)\} \leq D(G_1+G_2)\leq D(G_1)+D(G_2).$$
\end{theorem}
\proof
By contradiction, suppose that $D(G_1+G_2) < {\rm max}\{D(G_1),D(G_2)\}$. Without loss of generality we can assume that $D(G_1+G_2) < D(G_2)$. In this case the vertices of graph $G_2$ have been labeled with less than $D(G_2)$ labels, and so there exists a nontrivial automorphism $f_2$ of $G_2$ preserving the labeling of $G_2$. Hence there exists the nontrivial automorphism $h$ of $G_1+G_2$ preserving the labeling of $G_1+G_2$, which is contradiction.
\begin{equation*}
h(v) = \left\{
\begin{array}{ll} 
 v& \textsl{If} ~~v\in V(G_1),\\
f_2(v) &\textsl{If} ~~v\in V(G_2).
\end{array}\right.
\end{equation*}

To prove $D(G_1+G_2)\leq D(G_1)+D(G_2)$,  we first label $G_1$ in a distinguishing way with $D(G_1)$ labels, next we label the vertices of $G_2$ with the labels $\{D(G_1)+1,\ldots , D(G_1)+D(G_2)\}$ in a distinguishing way. This labeling is distinguishing because if  $f$ is an automorphism of $G_1+G_2$ preserving the labeling then with respect to the label of vertices of $G_1$ and $G_2$ we get that the restriction of automorphism $f$ to $G_i$ is $G_i$ where $i=1,2$, i.e.,  $f\vert_{G_1}=G_1$ and $f\vert_{G_2}=G_2$, and so $f\vert_{G_i}$ is an automorphism of $G_i$ for $i=1,2$. Since both  $G_1$ and $G_2$ have been labeled in a distinguishing way so we have $f\vert_{G_1}=id_{G_1}$ and $f\vert_{G_2}=id_{G_2}$. Therefore $f$ is the identity automorphism of $G_1+G_2$.\qed

\medskip 
To obtain a better upper bound for the distinguishing number of the join of two arbitrary graphs $G_1$ and $G_2$, we partition the vertices of $G_1+G_2$ such that every automorphism of $G_1+G_2$ maps the classes to each other. This partition is as follows:

\medskip
Let $G_1$ and $G_2$ be two graphs and $G=G_1+G_2$. Let $v_1$ be an arbitrary vertex of $G_1$. First put  $A_1=\overline{N_G(v_1)}$ (note that $\overline{N_G(v_1)} \subseteq V(G_1)$). 
We add all nonadjacent sets of the vertices of $G$ (say $v$) such that their nonadjacent sets satisfy $\overline{N_G(v)}\cap A_1 \neq \emptyset$, to $A_1$ and denote again the new  set by $A_1$ (if $v\in G$ and $\overline{N_G(v)}\cap A_1 \neq \emptyset$ then $v\in G_1$). We continue  this process  until there is no  vertex in $G$ with this property. 

Let $v_2$ be a vertex of $G_1$ such that $v_2\notin A_1$. Put    $A_2=\overline{N_G(v_2)}$ and similar to construction of $A_1$, add suitable nonadjacent sets of a vertex to $A_2$ and repeat this action. 
 It is clear that after a finite number of  steps, the vertices of $G_1$ partition to  $A_i$'s. With similar argument we suppose that  the vertices of $G_2$ partition to some sets, say, $B_j$'s. Without loss of generality we  assume that the vertices of $G$ are partitioned into $k+k'$ equivalence classes as follows (the notation $v$  is used for the vertices of $G_1$ and the notation $w$  is used for the vertices of $G_2$): 
\begin{align}\label{equation1}
&A_1=\overline{N_G(v_1)}\cup \cdots \cup \overline{N_G(v_{t_1})},\nonumber\\
&A_2=\overline{N_G(v_{t_1+1})}\cup \cdots \cup \overline{N_G(v_{t_1+t_2})},\nonumber\\
&\vdots \nonumber\\ 
&A_k=\overline{N_G(v_{t_1+\ldots + t_{k-1}+1})}\cup \cdots \cup \overline{N_G(v_{t_1+\ldots + t_{k}})},\\
&B_1=\overline{N_G(w_1)}\cup \cdots \cup \overline{N_G(w_{t'_1})},\nonumber\\
&B_2=\overline{N_G(w_{t'_1+1})}\cup \cdots \cup \overline{N_G(w_{t'_1+t'_2})},\nonumber\\
&\vdots \nonumber\\
&B_{k'}=\overline{N_G(w_{t'_1+\ldots + t'_{k'-1}+1})}\cup \ldots \cup \overline{N_G(w_{t'_1+\ldots + t'_{k'}})}.\nonumber
\end{align}

\begin{lemma}\label{lemma}
Let $G_1$ and $G_2$ be two graphs and $G=G_1+G_2$. Suppose that $ \mathcal{A}=\{A_1,\ldots , A_k\}$ and  $ \mathcal{B}=\{B_1,\ldots , B_{k'}\}$ are two partitions of the vertices $G_1$  and $G_2$ as stated in \eqref{equation1}, respectively. If $f$ is an automorphism of $G$, then $f$ is a permutation on the set $  \mathcal{A}\cup \mathcal{B}$.
\end{lemma}
\proof
 Let $u_1,u'_1\in V(G)$ and $f(u_1)=u'_1$. Since an automorphism preserves adjacency relation, $f(\overline{N_G(u_1)})=\overline{N_G(u'_1)}$.

Now  let $u_1,u_2,u'_1,u'_2\in V(G)$,  $f(\overline{N_G(u_1)})=\overline{N_G(u'_1)}$ and  $f(\overline{N_G(u_2)})=\overline{N_G(u'_2)}$. Then we have
\begin{equation*}
\overline{N_G(u_1)}\cap \overline{N_G(u_2)}\neq \emptyset\Leftrightarrow \overline{N_G(u'_1)}\cap \overline{N_G(u'_2)}\neq \emptyset .
\end{equation*}

 By induction, if $u_1,\ldots ,u_s,u'_1,\ldots ,u'_s\in V(G)$ and $f(\overline{N_G(u_i)})=\overline{N_G(u'_i)}$ where $1\leq i\leq s$ then we have
\begin{equation*}
\left(\overline{N_G(u_1)}\cup \ldots \cup \overline{N_G(u_s)}\right) \cap \overline{N_G(u_s)}\neq \emptyset\Leftrightarrow \left(\overline{N_G(u'_1)}\cup \ldots \cup \overline{N_G(u'_s)}\right) \cap \overline{N_G(u'_s)}\neq \emptyset .
\end{equation*}

By  the above illustrations and definitions of $A_i$ and $B_j$ with $1\leq i\leq k$ and $1\leq j\leq k'$ we can conclude that $f$ is a permutation on $  \mathcal{A}\cup \mathcal{B}$.\qed

\begin{corollary}\label{corj}
Let $G_1$ and $G_2$ be two graphs and $G=G_1+G_2$. Suppose that $ \mathcal{A}=\{A_1,\ldots , A_k\}$ and  $ \mathcal{B}=\{B_1,\ldots , B_{k'}\}$ are two partitions of the vertices $G_1$ and $G_2$ as stated in  \eqref{equation1}, respectively and put $   \mathcal{A}\cup  \mathcal{B}= \mathcal{C}=\{C_1,\ldots ,C_{k+k'}\}$. If $f$ is an automorphism of $G$ and $f(C_i)=C_j$, for some $i,j\in\{1,\ldots ,k+k'\}$, then induced subgraphs $G[C_i]$ and $G[C_j]$ are isomorphic. 
\end{corollary}

\begin{figure}
	\begin{center}
		\includegraphics[width=1\textwidth]{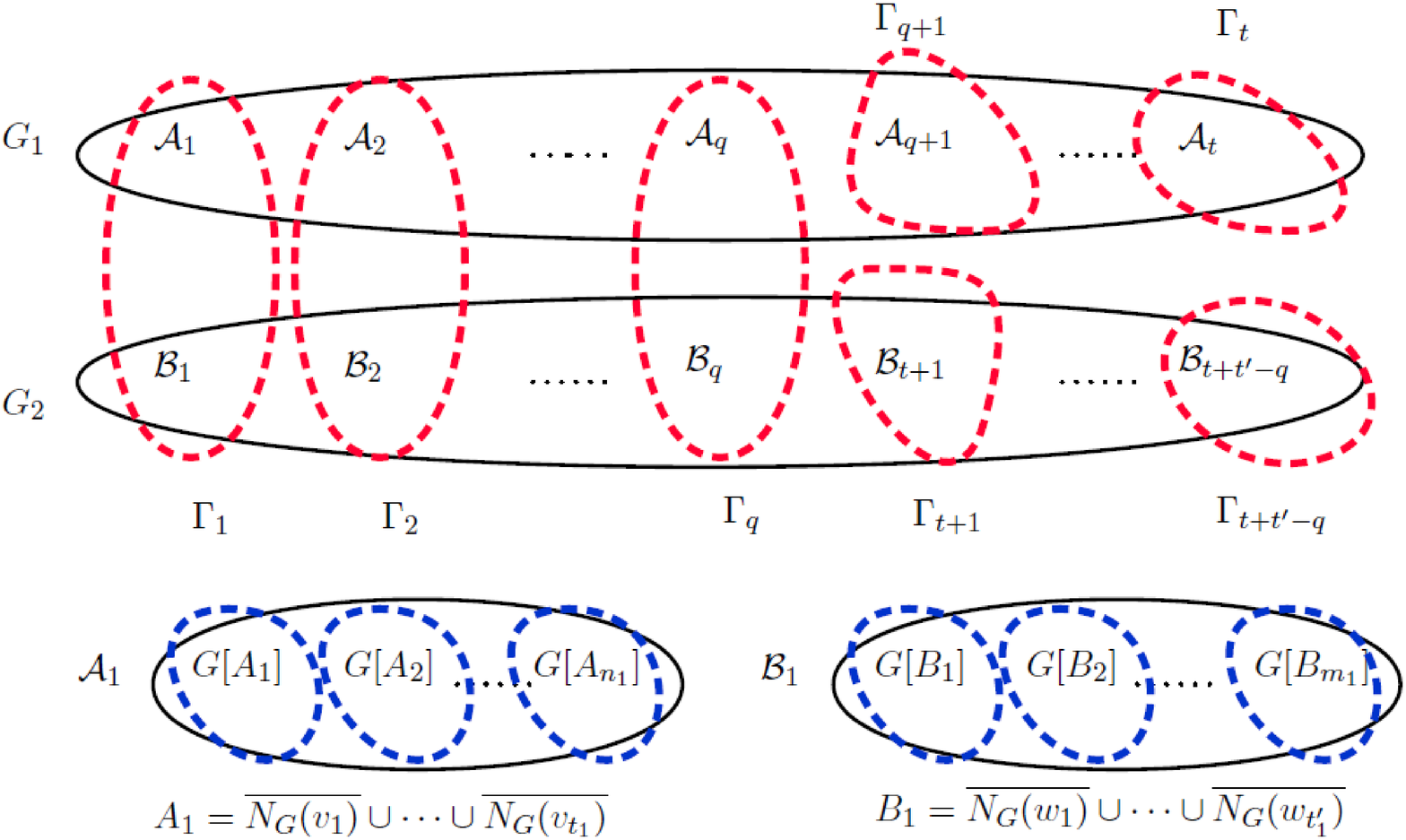}	
		\caption{The partition of $G_1+G_2$.}\label{figd}
	\end{center}
\end{figure}

Before stating and proving  the main theorems we need some additional information about $G_1$ and $G_2$.
Let $G_1$ and $G_2$ be two graphs and $G=G_1+G_2$ such that $ \mathcal{A}=\{A_1,\ldots , A_k\}$ and  $ \mathcal{B}=\{B_1,\ldots , B_{k'}\}$ are two partitions of the vertices $G_1$  and $G_2$ as stated in  \eqref{equation1}, respectively. Now we put  $H=\big\{G[A_1],\ldots ,G[A_k]\big\}$ and $H'=\big\{G[B_1],\ldots ,G[B_{k'}]\big\}$. Some of the induced subgraphs in each $H$ and $H'$ are isomorphic. We put all isomorphic induced subgraphs in $H$ and also $H'$, in a set and denote them by $\mathcal{A}_i$ and $\mathcal{B}_j$, respectively. In fact, we partitioned the two sets $H,H'$ into $t,t'$  disjoint sets $\mathcal{A}_1,\ldots , \mathcal{A}_t$ and $\mathcal{B}_1,\ldots , \mathcal{B}_{t'}$ such that $ \vert \mathcal{A}_i \vert = n_i $ and $ \vert \mathcal{B}_j \vert = m_j$ with $n_i,m_j\geq 1$, $1\leq i\leq t$ and $1\leq j\leq t'$ as follows:

\begin{align}\label{equation2}
&\mathcal{A}_1=\big\{G[A_1],\ldots , G[A_{n_1}]\big\},\nonumber\\
&\mathcal{A}_2=\big\{G[A_{n_1+1}],\ldots , G[A_{n_1+n_2}]\big\},\nonumber\\
&\vdots \nonumber\\ 
&\mathcal{A}_t=\big\{G[A_{n_1+\ldots +n_{t-1}+1}],\ldots , G[A_{n_1+\ldots +n_{t}}]\big\},\\
&\mathcal{B}_1=\big\{G[B_1],\ldots , G[B_{m_1}]\big\},\nonumber\\
&\mathcal{B}_2=\big\{G[B_{m_1+1}],\ldots , G[B_{m_1+m_2}]\big\},\nonumber\\
&\vdots \nonumber\\
&\mathcal{B}_{t'}=\big\{G[B_{m_1+\ldots +m_{t'-1}+1}],\ldots , G[B_{m_1+\ldots +m_{t'}}]\big\}.\nonumber
\end{align}
 It is possible that some of the elements in $\mathcal{A}_i$ are isomorphic to some elements in a $\mathcal{B}_j$, where $1\leq i \leq t$ and $1\leq j \leq t'$ (note that if an element of $\mathcal{A}_i$ is isomorphic to an element of $\mathcal{B}_j$ then all elements of $\mathcal{A}_i$ have this property). Let $q$ be  the number of $\mathcal{A}_i$ for which there exist some $\mathcal{B}_j$ that the elements of  $\mathcal{A}_i$ are isomorphic to elements of $ \mathcal{B}_j$. Then we can partition the set $H\cup H'$ into disjoint sets $\Gamma_1,\ldots, \Gamma_{t+t'-q}$ as follows: (we use new notation for vertices of $G$, if necessary).
 \begin{equation} \label{equation3}
\left\{\begin{array}{ll}
\Gamma_i=\mathcal{A}_i\cup \mathcal{B}_i&~1\leq i\leq q,\\
\Gamma_{q+i}=\mathcal{A}_{q+i}&~1\leq i\leq t-q,\\
\Gamma_{t+i}=\mathcal{B}_{t+i}&~1\leq i\leq t'-q,\\
\end{array}\right.
\end{equation}
where $0\leq q\leq {\rm min}\{t,t'\}$ (See Figure \ref{figd}).

\begin{remark}\label{rem}
Using the partition of $H\cup H'$ in \eqref{equation3}, Lemma \ref{lemma} and Corollary \ref{corj} we can conclude that  if $f\in {\rm Aut}(G_1+G_2)$ then $f\vert_{\Gamma_i} \in {\rm Aut}(\Gamma_i)$ for $1\leq i\leq q$ where $0\leq q\leq {\rm min}\{t,t'\}$.
\end{remark}

Now we are ready to state and prove the main result on the distinguishing number of the join of two graphs: 
\begin{theorem}\label{Djoin}
Let $G_1$  and $G_2$ be two non-isomorphic graphs and $G=G_1+G_2$.
\begin{enumerate}
\item[(i)] If $q=0$ then $D(G_1+G_2)={\rm max}\{D(G_1),D(G_2)\}$.
\item[(ii)] If $q\neq 0$ and $z={\rm min}\big\{{\rm max}\{n_1,\ldots , n_q\}, {\rm max}\{m_1,\ldots ,m_q\}\big\}$ then $$D(G_1+G_2)\leq {\rm max}\{D(G_1),D(G_2)\}+z.$$
\end{enumerate}
\end{theorem}
\proof
\begin{enumerate}
	\item[(i)] 
	 If $q=0$, then there is no element of $H$  isomorphic to an element of $H'$. By Corollary \ref{corj}, if $f\in {\rm Aut}(G)$ then $f\vert_{G_1}\in {\rm Aut}(G_1)$ and  $f\vert_{G_2}\in {\rm Aut}(G_2)$, and so $D(G_1+G_2)\leq {\rm max}\{D(G_1),D(G_2)\}$. Therefore by Theorem \ref{thh5} we have the result.

\item[(ii)] Let $d={\rm max}\{D(G_1),D(G_2)\}$. We shall present a distinguishing labeling with $d+z$ labels. Without loss of generality we can assume that $z=m_1$, and so $\mathcal{B}_1=\{G[B_1],\ldots ,G[B_{m_1}]\}$.

First, we label both $G_1$ and $G_2$ with $D(G_1)$ and $D(G_2)$ labels in a distinguishing way, respectively.  Now to obtain a distinguishing labeling of $G_1+G_2$, we change the labels of the vertices $G_2$ as follows:
 \begin{itemize}
 \item We  change the label of an arbitrary vertex of  $G[B_i]$ to $d+i$,  for every $1\leq i\leq m_1$.
 \end{itemize}
 We do similar above process on $\mathcal{B}_2,\ldots , \mathcal{B}_{t'}$ (note that if $z=n_{i_k}$ for some $k\in \{1,\ldots , q\}$ then we should do the similar work  on $G_1$).  By Lemma \ref{lemma}, Corollary \ref{corj} and the distinguishing labeling in both $G_1$ and $G_2$,  we can conclude that presented labeling is distinguishing. Since we used ${\rm max}\{D(G_1),D(G_2)\}+z$ labels, the inequality follows.\qed
 \end{enumerate} 

\begin{remark}
The value of $z$ in Theorem \ref{Djoin} (ii) can be zero or  sufficiently large,  depending on the structure of graphs $G_1$ and $G_2$. As an example, consider the  complete
 $k$-partite  graph  $K_{d, \ldots ,d}$ as $G_1$ and $G_2$ and  $G= K_{d, \ldots ,d}+K_{d, \ldots ,d}$, then using notations in \eqref{equation2},  $\mathcal{A}_i=\mathcal{B}_i=\emptyset$  for $2\leq i\leq t$, and  $\mathcal{A}_1=\mathcal{B}_1=\{G[A_1],\ldots , G[A_k]\}$, where $A_i$ is the $i$-th part of $K_{d, \ldots ,d}$. Therefore $z = k$ and so, $z$ can be sufficiently large.
\end{remark}

Now we shall  show that the inequality in Theorem \ref{Djoin} (ii) is sharp.
\begin{corollary}\label{exa}
Let $n>m,n>m'$ and $m\neq m'$. The distinguishing number of  $K_{n,m}+K_{n,m'}$ is $n+1$.
\end{corollary} 
\proof
  Let $X=\{v_1,\ldots , v_n\}$,  $Y=\{w_1,\ldots , w_m\}$ be two parts of $K_{n,m}$, and $X'=\{v'_1,\ldots , v'_n\}$, $Y'=\{w'_1,\ldots , w'_{m'}\}$ be two parts of $K_{n,m'}$. Suppose that  $G= K_{n,m}+K_{n,m'}$. Using the partition in \eqref{equation1} we can write:
\begin{equation*}
A_1=\overline{N_G(v_1)}=\{v_1,\ldots , v_n\},~~~A_2=\overline{N_G(w_1)}=\{w_1,\ldots , w_m\}.
\end{equation*}
 Since the number of elements in $A_1$ and $A_2$ are distinct, $G[A_1]\ncong G[A_2]$. Then by the partition in \eqref{equation2} we have 
  $\mathcal{A}_1=\{G[A_1]\}$ and $\mathcal{A}_2=\{G[A_2]\}$, and so $n_1=n_2=1$. Now by similar argument we can write:
\begin{equation*}
B_1=\overline{N_G(v'_1)}=\{v'_1,\ldots , v'_n\},~~~B_2=\overline{N_G(w'_1)}=\{w'_1,\ldots , w'_{m'}\}.
\end{equation*}
Then $\mathcal{B}_1=\{G[B_1]\}$ and $\mathcal{B}_2=\{G[B_2]\}$, and so $m_1=m_2=1$.  Since the induced subgraphs have no edges, $G[A_1]\cong G[B_1]$. With respect to the partition in \eqref{equation3} we have 
\begin{equation*}
\Gamma_1=\mathcal{A}_1\cup \mathcal{B}_1=\{G[A_1],G[B_1]\},~\Gamma_2= \mathcal{A}_2=\{G[A_2]\},~\Gamma_3=\mathcal{B}_2=\{G[B_2]\}.
\end{equation*}
  It is clear that for every labeling by $n$ labels we can find a labeling preserving automorphism of $\Gamma_1$.   So we can find an automorphism of $G$ with this property. Consider the following labeling by $n+1$ labels:
 
 We assign to the vertices in $A_1$ the labels $1,\ldots ,n$ and to the vertices in $B_1$ the labels $1,\ldots ,n-1,n+1$. We label the vertices in $A_2$ with the labels $1,\ldots ,m$ and the vertices in $B_2$ with the labels $1,\ldots ,m'$. By Remark \ref{rem}, this labeling is distinguishing, and so 
 $D( K_{n,m}+K_{n,m'})=n+1$. \qed
 
\begin{theorem}\label{thh2.11}
 Let $n_1,...,n_t$ be the number of elements of classes stated in \eqref{equation2}. We have 
  $$ D(G)\leq D(G+G)\leq D(G)+{\rm max}\{n_1,\ldots ,n_t\}.$$
\end{theorem}
\proof
Let $G_1$ and $G_2$ be two isomorphic graphs and denote both of them  by $G$, then the left side inequality is identified by Theorem \ref{thh5}. To   prove  the right side of inequality, we present a distinguishing labeling as follows: 
 
 Without loss of generality we can assume that $n_1={\rm max}\{n_1,\ldots ,n_t\}$. First, we label $G$ and its copy with $D(G)$ labels in a distinguishing way.  To obtain a distinguishing labeling for $G+G$ we change the labels of the vertices of $G$ as follows:
 \begin{itemize}
 \item We  change the label of an arbitrary vertex of  $(G+G)[A_i]$ to $D(G)+i$,  for every $1\leq i\leq n_1$.
 \end{itemize}
So the labels of vertices of $\mathcal{A}_1$ were changed. We do similar  process on $\mathcal{A}_2,\ldots , \mathcal{A}_t$.  By Lemma \ref{lemma}, Corollary \ref{corj} and the distinguishing labeling in both $G$ and its copy,  we can conclude that presented labeling is distinguishing. Since we used $D(G)+{\rm max}\{n_1,\ldots , n_t\}$ labels, the right side inequality follows.\qed

\begin{remark} 
	With similar argument as in the proof of  Corollary \ref{exa} we can show that  the inequality in Theorem \ref{thh2.11} is sharp for the star graphs $K_{1,n}$. In fact $D(K_{1,n}+K_{1,n})=n+1$ where $  D(K_{1,n})=n$ and $ {\rm max}\{n_1,\ldots , n_t\}=1$.
\end{remark}

\section{Distinguishing index of the join of two graphs} 

In this section  we study the distinguishing index for the join of two graphs. We say that a graph $G$ is almost spanned by a subgraph H if $G-v$ is spanned
by $H$ for some $v \in V (G)$. We need the following lemmas 
in this section.

\begin{lemma}{\rm \cite{nord}}\label{nordspannin}
 If a graph $G$ is spanned or almost spanned by a subgraph $H$, then $D'(G) \leq D(H) + 1$.
\end{lemma}

\begin{lemma}{\rm \cite{nord}}\label{traceable}
Let $G$ be a graph of order $n \geq 7$ with a Hamiltonian path, then $D'(G) \leq 2$.
\end{lemma}

By these two lemmas, we can obtain the following upper bounds for the distinguishing index of join of two graphs.

\begin{theorem}
Let  $G$ and $H$ be two graphs of orders  $n$ and $m$, respectively. Then $D'(G + H) \leq D'(K_{n,m})+1$.
\end{theorem}
\proof Since  the  complete bipartite graph  $K_{n,m}$, is  a spanning subgraph $G + H$, we can conclude the result by Lemma \ref{nordspannin}. \qed

\begin{theorem}\label{disindrestorder}
If $G$ has $n$ vertices and $H$ has $m$ vertices, such that $4 \leq n \leq m \leq 2n$, then $D'(G+H)\leq 2$.
\end{theorem}
\proof We  use the complete bipartite $K_{n,m}$ subgraph to find an asymmetric
spanning subgraph of $G + H$. Now we have the result by Lemma \ref{nordspannin}. \qed

\begin{theorem}\label{indjoindelta}
Let $G$ and $H$ be two graphs of orders $n$ and $m$, respectively such that   $\delta(G)\leq \delta(H)$.  If ${\rm {\rm min}}\{\delta(G) + m, \delta(H)+n\}\geq \frac{n+m-1}{2}$ and $m+n\geq 7$, then $D'(G+H)\leq 2$.
\end{theorem}
\proof It is known that if the {\rm min}imum degree of a graph of order $n$ is at least $\frac{n-1}{2}$, then graph has a Hamiltonian path. Since the {\rm min}imum degree of $G+H$ is  ${\rm {\rm min}}\{\delta(G) + m, \delta(H)+n\}$, so the result follows by Theorem \ref{traceable}.\qed

\begin{corollary}
If $G$ is a graph of order $n\geq 2$, then $D'(\underbrace{G+ \cdots +G}_{k-times})=2$ for any $k \geq 2$, except $D'(K_2 +K_2) =3$.
\end{corollary}
\proof  For $k =2$,  we have $\delta (G+G)= \delta (G)+n\geq \frac{2n-1}{2}=\frac{|G+G|-1}{2}$, and hence $G+G$ has a Hamiltonian path.  If $n \geq 4$, then  $2n=|G+G|\geq 7$,  and so we  have  $D'(G+G)\leq 2$, by Lemma \ref{traceable}. On the other hand, since the automorphism group of graph $G+G$ is non-trivial, so $D'(G+G)\geq 2$.  Therefore $D'(G+G)= 2$. If $n=3$, then it is easy to see that $D'(G+G)=2$. 
Now a simple induction argument together with Theorem \ref{indjoindelta} yield that $D'(\underbrace{G+ \cdots +G}_{k-times})=2$, for any $k \geq 2$. \qed

To obtain an upper bound for $D'(G_1+G_2)$ we consider \eqref{equation3} which is a partition of $H\cup H'$, i.e., $\Gamma_1,\ldots , \Gamma_{t+t'-q}$. Note that the elements of $\mathcal{A}_i$ are isomorphic to elements of  $\mathcal{B}_i$ for $1\leq i\leq q$ where $0\leq q\leq {\rm min}\{t,t'\}$. If $G_1\cong G_2$ then $t=t'$ and the elements of $\mathcal{A}_i$ are isomorphic to elements of  $\mathcal{B}_i$ for $1\leq i\leq t$.

Let $E_i$ be the set of edges of $G_1+G_2$ such that the end points of its edges are in $\Gamma_i$ for  $1\leq i\leq t+t'-q$. We add the set $E_i$ to the set of edges $\Gamma_i$ and denote  the  obtained new graph by  $\Gamma'_i$. The following result gives an upper bound for $D'(G_1+G_2)$ based on the distinguishing index of $\Gamma'_i$. 

\begin{theorem}\label{thmd'1}
Let $G_1$ and $G_2$ be two graphs such that $G_1+G_2$ has been partitioned to the set of induced subgraphs $\Gamma_1,\ldots , \Gamma_{t+t'-q}$ as \eqref{equation3}. Then $$D'(G_1+G_2)\leq {\rm max}\{D'(\Gamma'_1),\ldots , D'(\Gamma'_{t+t'-q})\}.$$ 
\end{theorem}
\proof We label the edges of the graph $\Gamma'_i$ ($1\leq i\leq t+t'-q$)  by $D'(\Gamma'_i)$ labels in a distinguishing way. We assign the remaining edges the label $1$. By Remark \ref{rem}, this labeling is distinguishing. The number of labels that have been used here is
\begin{equation*}
{\rm max}\{D'(\Gamma'_1),\ldots , D'(\Gamma'_{t+t'-q})\}.
\end{equation*}
So we have the result.\qed

Now, we like to present another upper bound for $D'(G_1+G_2)$. For this purpose we state some preli{\rm min}aries. 

\medskip 

Let $X_i$, $i\in I$ ($I$ is the index set) be the set of complete bipartite graphs $K_{|V(\Gamma_s)|,|V(\Gamma_{s'})|}$  satisfying   following two conditions:
\begin{itemize}
\item[$\bullet$] The two parts of each element of $X_i$ should be distinct.
\item[$\bullet$] The set of all parts that have been used as parts of  elements of $X_i$ should be $\{V(\Gamma_1),\ldots , V(\Gamma_{t+t'-q})\}$. 
\end{itemize}

Let $\varepsilon_i={\rm max}\{D'(K_{|V(\Gamma_s)|,|V(\Gamma_{s'})|}):K_{|V(\Gamma_s)|,|V(\Gamma_{s'})|}\in X_i\}$. Then we have the following theorem:

\begin{theorem}\label{thmd'2}
Let $G_1$ and $G_2$ be two graphs such that $G_1+G_2$ has been partitioned to induced subgraphs  $\Gamma_1,\ldots , \Gamma_{t+t'-q}$ as \eqref{equation3}. Then $D'(G_1+G_2)\leq {\rm min}\{\varepsilon_i\}_{i\in I}$. 
\end{theorem}
\proof
 We label the edges of each complete bipartite graph in $X_i$ in  distinguishing way (by $D'(K_{|V(\Gamma_s)|,|V(\Gamma_{s'})|})$ labels) and assign to the remaining edges the label $1$. Since all parts $\Gamma_1,\ldots , \Gamma_{t+t'-q}$ have been used in building of the complete bipartite graphs in $X_i$ and by Remark \ref{rem}, this labeling is distinguishing. Therefore $D'(G_1+G_2)\leq {\rm min}\{\varepsilon_i\}_{i\in I}$.\qed

%The case $ \gamma =0$ has been considered in Theorem \ref{thh2}.

\begin{remark}\label{compare} 
By setting $\lambda_1={\rm max}\{D'(\Gamma'_1),\ldots , D'(\Gamma'_{t+t'-q})\}$ and $\lambda_2={\rm min}\{\varepsilon_i\}_{i\in I}$ and by Theorem \ref{thmd'1} and \ref{thmd'2} we have  $D'(G_1+G_2)\leq {\rm min}\{\lambda_1,\lambda_2\}$. This  raises the question ``which upper bound is better, $ \lambda_1$ or $\lambda_2$"? We show that for some graphs the upper bound $\lambda_1$ is better than $\lambda_2$ and for some graphs the situation is different. We present two examples and these    examples show also that the upper bounds of the Theorem  \ref{thmd'1} and \ref{thmd'2} are sharp.
\end{remark}

Since the line graph of $K_{k,n}$ is isomorphic to Cartesian product  $K_k\Box K_n$, so  ${\rm Aut}(K_{k,n})$ coincides with ${\rm Aut}(K_k\Box K_n)$. Therefore  the distinguishing index of the complete bipartite graphs which is needed in the solution of Example \ref{pathjoin} can be translated to distinguishing number of Cartesian product of complete graphs. 
 
 \begin{theorem}\label{thmimrich}{\rm\cite{Imrich}}
 Let $k,n,d$ be integers so that $d\geq 2$ and $(d-1)^k< n\leq d^k$. Then
 \begin{equation*}
D(K_k\Box K_n)=\left\{
 \begin{array}{ll}
 d&\textsl{If $n\leq d^k -\lceil {\rm log}_d k\rceil -1$,}\\
 d+1&\textsl{If $n\geq d^k -\lceil {\rm log}_d k\rceil +1$.}
 \end{array}\right.
 \end{equation*}
 If $n=d^k -\lceil {\rm log}_d k\rceil$ then $D(K_k\Box K_n)$ is either $d$ or $d+1$ and can be computed recursively in $O({\rm log}^*(n))$ time.
 \end{theorem}

\begin{example}\label{pathjoin} 
The upper bound  in Theorem \ref{thmd'1} is better than the upper bound  in Theorem \ref{thmd'2} for the $D'(P_n+P_m)$ with $n,m\geq 2$ and $ n\neq m$. 
\end{example}
{\bf Solution.}  
Set $G=P_n+P_m$. Suppose that $V(P_n)=\{v_1,\ldots , v_n\}$ and $V(P_m)=\{w_1,\ldots , w_m\}$. With these notations  we have $A_{1}=\overline{N_{G}(v_1)}=\{v_1, \ldots , v_n\}$ and $B_{1}=\overline{N_{G}(w_1)}=\{w_1, \ldots , w_m\}$. Thus  $\mathcal{A}_1=\{G[A_1]\}$ and $\mathcal{B}_1=\{G[B_1]\}$. Since $ n\neq m$, so $\Gamma_1= \mathcal{A}_1=\{G[A_1]\}$,  $\Gamma_2=\mathcal{B}_1=\{G[B_1]\}$ and $q=0$. Also, $ \Gamma'_1=P_n$ and $ \Gamma'_2=P_m$. If we label both  $ \Gamma'_1$ and $ \Gamma'_2$ by two labels in a  distinguishing way (note that $D'(P_n)=D'(P_m)=2$) then we have a distinguishing labeling with two labels by Remark \ref{rem}. 

It is easy to see that  $D'(P_n+P_m)=\lambda_1=2$, and so the inequality of Theorem \ref{thmd'1} is sharp. On the other hand, using the notation of Theorem \ref{thmd'2} we have $I=\{1\}$, and so $X_1=\{K_{|V(\Gamma_1)|,|V(\Gamma_2)|}\}$. By Theorem \ref{thmimrich} it is clear that $\varepsilon_1=D'(K_{V(\Gamma_1),V(\Gamma_2)})$ is not equal with $2$ for all $m,n\geq 2$. Therefore the upper bound $\lambda_1$ is better than $\lambda_2$ for $D'(P_n+P_m)$.\qed
\medskip

\begin{figure}
	\begin{center}
	\includegraphics[width=0.45\textwidth ]{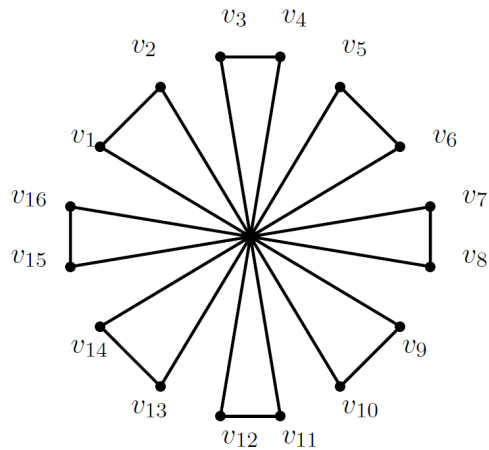}	
\caption{The graph  $F_8$.}\label{friend}
	\end{center}
\end{figure}

Here we shall present two graphs for which the upper bound in Theorem \ref{thmd'2} is better that the upper bound in Theorem \ref{thmd'1}.
We recall that the friendship graph $F_n$ is the join of $K_1$ with $nK_2$. In other words, $F_n$ can be constructed by joining $n$ copies of the cycle graph $C_3$ with a common vertex  (see Figure \ref{friend}). The following theorem gives  the distinguishing index of the friendship graph $F_n$.

\begin{theorem}{\rm\cite{soltani}} 
Let $a_n=1+27n+3\sqrt{81n^2+6n}$. The distinguishing index  of the friendship graph $F_n$  $(n\geq 2)$ is
 $$D'(F_n)=\lceil\frac{1}{3} (a_n)^{\frac{1}{3}}+\frac{1}{3(a_n)^{\frac{1}{3}}}+\frac{1}{3}\rceil.$$
	\end{theorem}

\begin{example}
The upper bound in Theorem \ref{thmd'2} is better than the upper bound in Theorem \ref{thmd'1} for   $D'(F_n+F_m)$, where $ 2\leq n < m$. 
\end{example} 
{\bf Solution.}	
 Suppose that  $G=F_n+F_m$. The central vertices of $F_n$ and $F_m$ are denoted by $x_0$ and $y_0$, respectively. Any two adjacent vertices of $F_n$  (except central vertex $x_0$) are denoted by $x_{2i-1}$ and $x_{2i}$ where $i=1,\ldots ,n$. The corresponding vertices of $F_m$ are denoted by $y_{2j-1}$ and $y_{2j}$ where $j=1,\ldots ,m$.

By the partition in \eqref{equation1} we can write  $A_1=\overline{N_G(x_0)}=\{x_0\}$ and $A_2=\overline{N_G(x_1)}\cup \overline{N_G(x_2)} =\{x_1,\ldots ,x_{2n}\}$. Also $B_1=\overline{N_G(y_0)}=\{y_0\}$ and $B_2=\overline{N_G(y_1)}\cup \overline{N_G(y_2)}=\{y_1,\ldots ,y_{2m}\}$. By the partition in \eqref{equation2}, $\mathcal{A}_1=\{G[A_1]\}$ and $\mathcal{A}_2=\{G[A_2]\}$, also $\mathcal{B}_1=\{G[B_1]\}$ and $\mathcal{B}_2=\{G[B_2]\}$. Let $m\neq n$.  By \eqref{equation3} and the hypothesis $m\neq n$ we have, $\Gamma_1=\mathcal{A}_1\cup \mathcal{B}_1=\{G[A_1],G[B_1]\}$,  $\Gamma_2=\mathcal{A}_2=\{G[A_2]\}$ and $\Gamma_3=\mathcal{B}_2=\{G[B_2]\}$, and so $q=1$. By notation of Theorem \ref{thmd'2}, one of the sets $X_i$ is  $X_1=\{K_{|V(\Gamma_1)|,|V(\Gamma_2)|}, K_{|V(\Gamma_2)|,|V(\Gamma_3)|}\}$. 
By Theorem \ref{thmimrich}, we have 
\begin{align*}
\varepsilon_{1}&={\rm max}\{D'({K_{|V(\Gamma_1)|,|V(\Gamma_2)|}}),D'(K_{|V(\Gamma_2)|,|V(\Gamma_3)|})\}\\
&= {\rm max}\{D'({K_{2,2n}}),D'(K_{2n,2m})\}.
\end{align*}

 Thus $\lambda_2\leq  \varepsilon_{1}$.  On the other hand, $\Gamma'_i$, $i\in \{2,3\}$  is the union of graphs $P_2$, and so the distinguishing index of graphs $\Gamma'_2$ and $\Gamma'_3$ has not defined. Therefore the upper bound $\lambda_2$ is better than $\lambda_1$.\qed

%%%%%%%%%%%%%%%%%%%%%%%%%%%%%%%%%%%%%%%%%%%


\begin{thebibliography}{99}
 	\bibitem{Albert2005} M.O. Albertson, {\it Distinguishing Cartesian product of graphs}, Electron. J. Combin. 12 (2005) \#N17. 
 	
 	\bibitem{Albert} M.O. Albertson and K.L. Collins, {\it Symmetry breaking in graphs}, Electron. J. Combin. 3 (1996) \#R18.
 	
 	\bibitem{soltani} S. Alikhani and S. Soltani, {\it Distinguishing number and distinguishing index of certain graphs}, Filomat, to appear. Available at \texttt{http://arxiv.org/abs/1602.03302}.
 	 
 	\bibitem{bogs} 
 	B. Bogstad, L. Cowen, J. Lenore,  {\it The distinguishing number of the hypercube},
 	Discrete Math. 283 (1) (2004) 29-35. 
 	 
 	
 	\bibitem{Sandi} R. Hammack, W. Imrich and S. Kalav\v zar, {\it Handbook of product graphs (second edition)}, Taylor \& Francis group (2011).  
 
 
  	\bibitem{Imrich} W. Imrich, J. Jerebic and S. Klav\v zar, {\it The distinguishing number of Cartesian products of complete graphs}, European J. Combin. 29 (4), (2008) 922-929.
 	
 	\bibitem{Kali1} R. Kalinowski and M. Pilsniak, {\it Distinguishing graphs by edge colourings}, European J. Combin. 45(2015) 124-131.
 	
 	\bibitem{Klavzar} S. Klav\v zar and X. Zhu, {\it Cartesian powers of graphs can be distinguished by two labels}, European J. Combin. 28 (2007) 303-310.  
 
 
  	\bibitem{fish} F. Michael and I. Garth, {\it Distinguishing colorings of Cartesian products of complete graphs}, Discrete Math., 308 (11), (2008) 2240-2246. 
 	
 	\bibitem{nord}  M. Pil\'sniak, {\it Nordhaus-Gaddum bounds for the distinguishing index}, Available at \texttt{www.ii.uj.edu.pl/preMD/}. 
 	

 	
 	

 	
 	
 \end{thebibliography}
\end{document}